\documentclass[12pt]{article}

\usepackage{amssymb}

\usepackage[utf8]{inputenc}
\usepackage[T1]{fontenc}
\usepackage{pgf,pgfarrows,pgfnodes,pgfautomata,pgfheaps}
\usepackage{amscd}
\usepackage{amsfonts}
\usepackage{latexsym}
\usepackage{graphicx}
\usepackage{amsthm}
\usepackage{amsmath}
\usepackage{leftidx}

\theoremstyle{definition}
\newtheorem{thm}{Theorem}

\begin{document}

\date{}
\author{}

\begin{center}
\textbf{\Large{On the differential transcendentality of the Morita p-adic gamma function }}\\
\vspace{10pt}
EL\.ZBIETA ADAMUS \footnote{This work was partially supported by the Faculty of Applied Mathematics AGH UST
statutory tasks within subsidy of Ministry of Science and Higher Education.} \\ Faculty of Applied Mathematics, \\ AGH University of Science and Technology \\
al. Mickiewicza 30, 30-059 Krak\'ow, Poland \\
e-mail: esowa@agh.edu.pl 

\end{center}

\vspace{1cm}
\begin{abstract}
 In this note we prove that for any given prime number $p$ the Morita $p$-adic gamma function $\Gamma_p$ is differentially transcendental over $\mathbb{C}_p(X)$.
\end{abstract}

\vspace{15pt}

\section{Introduction}

\par For any differential field $K$ denote by $K\{Y\}$ the ring of differential polynomials with coefficients in $K$.
If $K \subset L$ is a differential field extension, then any element $f \in L$ is called \emph{differentially algebraic} over $K$ if and only if there is a non-zero differential polynomial $P \in K\{Y\}$ such that $P(f) = P(f,f',f'', \ldots, f^{(m)})=0$.  If $f$ is not differentially algebraic over $K$ then we say that $f$ is \emph{differentially transcendental} over $K$.

One can consider the field of complex-valued rational functions $\mathbb{C}(z)$ of single complex variable $z$ as a differential field equipped with the standard derivation $\frac{d}{dz}$. 
Functions which are   differentially transcendental over $ \mathbb{C}(z) $, are also called
\emph{transcendentally transcendental}  (term introduced by E. H. Moore in 1896) or \emph{hypertranscendental} (term introduced by D. D. Morduhai-Boltovskoi in 1914).
A well known example of a differentially transcendental function is \textit{Euler's gamma function} $\Gamma(z)$. If we consider the \emph{Euler integral of the second kind}, i.e. an improper integral of the form \[\int_0^{\infty} t^{z-1}e^{-t}dt,\]
 it converges absolutely on $\{z \in \mathbb{C}: Re{(z)} >0\}$.
 For any $z \in \mathbb{C}$ such that $Re{(z)} >0$  we define
\[ \Gamma(z) = \int_0^{\infty} t^{z-1}e^{-t}dt. \] 
Using the equality $\Gamma(z) = \frac{\Gamma(z+1)}{z}$ one can extend the definition on the set $\mathbb{C} \setminus \mathbb{Z}_{\leq 0}$, where $\mathbb{Z}_{\leq 0}=\{k \in \mathbb{Z}, k \leq 0\}$. Indeed we have
$\Gamma(z+n)=(z-n+1)\ldots(z+1)z \cdot \Gamma (z)$. 
Observe that $\Gamma(z+n)$ is holomorphic for $z \in \mathbb{C}$ such that $Re(z)>-n$.
We use the identity theorem for holomorphic functions 
to obtain the meromorphic function $\Gamma : \mathbb{C} \setminus \mathbb{Z}_{\leq 0} \rightarrow \mathbb{R}_{+}$. Hölder's theorem states that $\Gamma$ does not satisfy any algebraic differential equation whose coefficients are rational functions. For the proof see for example \cite{T}.

\begin{thm}[Hölder, 1887] 
 The gamma function $\Gamma$ is differentially transcendental over $\mathbb{C}(z)$.
\end{thm}

\section{The Morita p-adic gamma function}

\par In $p$-adic analysis one can find an analog of the classical gamma function defined on the ring of $p$-adic integers. In 1975 in \cite{Mo} Morita defined the p-adic  gamma function $\Gamma _p$ by a suitable modification of
the function $n \mapsto n!$.  
 
We denote by $ \mathbb{Z}_p, \mathbb{Q}_p$ the ring of $p$-adic integers and the field of $p$-adic numbers, respectively. By $\mathbb{Z}_p^{\times}$ we denote the group of invertible elements in $\mathbb{Z}_p$
and by $\mathbb{C}_p$ we denote the completion of the algebraic clousure of the field $\mathbb{Q}_p$.
The factorail function $n \mapsto n!$ cannot be extended by continuity to a function $f: \mathbb{Z}_p \rightarrow \mathbb{Q}_p$ such that $f(n)=n f(n-1)$ for all integers greater than $0$. When $p$ is odd, one need to consider the \emph{restricted factorial} $n!*$ given by
\[ n!*=\prod_{\substack{1 \leq j \leq n \\ p \, \nmid \, j }}j\]
and use the generalization of the \emph{Wilson congruence}
$ (p-1)! \equiv -1 \mod p$.
\emph{The Morita $p$-adic gamma function} is the continuous function $f: \mathbb{Z}_p \rightarrow \mathbb{Z}_p$ extending $f(n):=(-1)^n n!*$, for $n \geq 2$. One can observe that the values of $\Gamma_p$ belong to $\mathbb{Z}_p^{\times}$.  Recall that we have the partition 
$\mathbb{Z}_p= p\mathbb{Z}_p \cup \mathbb{Z}_p^{\times}$, where  $p\mathbb{Z}_p=\{x \in \mathbb{Z}_p: \, |x|_p<1\}$ is the only maximal ideal  of the ring $\mathbb{Z}_p$. What is more 
\begin{equation}
 \Gamma_p(x+1)=h_p(x) \cdot \Gamma_p(x), \quad \mathrm{where} \quad h_p(x) = \left\{ \begin{array}{lll}
                                                                                 -x&;&x \in \mathbb{Z}_p^{\times}\\
                                                                                 -1&;& x \in p \mathbb{Z}_p
                                                                                \end{array} \right. .
 \label{gam1}
\end{equation}
It can be proved that all above formulas hold also for $p=2$. For more detailed information see for example \cite{R}, section 7.1.

\section{Main result}

Consider the field $\mathbb{C}_p(X)$ of rational functions with coefficients in $\mathbb{C}_p$.
We claim the following.

\begin{thm} For a given prime number $p$, the
 Morita $p$-adic gamma function $ \Gamma_p$ is differentially transcendental over $\mathbb{C}_p(X)$.
 \label{CPH}
\end{thm}

 In the following proof we consider the \emph{antilexicographical ordering}
of monomials in $\mathbb{C}_p(X)[Y_0,Y_1, \ldots, Y_n]$. If we consider $\alpha, \beta \in \mathbb{Z}^{n+1}$, $\alpha \geq 0, \beta \geq 0$, then 
$\alpha > \beta$ if in the vector difference $\alpha - \beta$ the leftmost nonzero entry is negative. For 
$Y=(Y_0,Y_1, \ldots, Y_n)$ we write $Y^{\alpha} > Y^{\beta}$ if $\alpha > \beta$. Thus we have $Y_0 < Y_1 < \ldots < Y_n$.\\

\textit{Proof of theorem \ref{CPH}.} Suppose on the contrary, that  $\Gamma_p$ is not differentially transcendental over $\mathbb{C}_p(X)$. 
If $\Gamma_p$ is differentially algebraic over $\mathbb{C}_p(X)$, then so is its restriction $g:=\Gamma_p|_{p\mathbb{Z}_p} : p\mathbb{Z}_p\rightarrow \mathbb{Z}_p^{\times}$. Hence there exists a differential polynomial $P \in \mathbb{C}_p[X]\{Y\}$ such that $P=P(X,Y_0,Y_1, \ldots, Y_n) \not\equiv 0$ and 
\begin{equation}
 \forall x \in p \mathbb{Z}_p \quad P(x, g(x),  g'(x), \ldots,  g^{(n)}(x))=0.
 \label{eqP}
\end{equation}
Without loss of generality, one can assume that $P$ contains a monomial having a non-zero power of one of the indeterminates $Y_0, Y_1,\ldots, Y_n $. This is due to the fact that  $P$ cannot define a function which is zero on $g (p \mathbb{Z}_p)$. 
Let $q(X)Y_0^{a_0}Y_1^{a_1}\ldots Y_n^{a_n}$, where $q \in \mathbb{C}_p(X)$,  be the leading term of $P$, i.e. the one with the biggest $(a_0,a_1, \ldots, a_n  ) \in \mathbb{Z}^{n+1}$ with respect to the antilexicographical ordering. Let us denote it by $LT(P)$. Assume that from all such polynomials we choose $P$, which is minimal in the sense that it has a minimal leading term $LT(P)$
 and moreover $\deg q$ is minimal and the leading coefficient of $q$ is $1$.
 
For every $x \in p \mathbb{Z}_p$  we have $1= |1|_p \neq|x|_p <1$. Hence $|x+1|_p = \max\{|x|_p, 1\}=1$, which is equivalent to
$x+1 \in \mathbb{Z}_p^{\times}$. By (\ref{gam1}) for every $x \in p\mathbb{Z}_p $ we obtain
\[
 \Gamma_p(x+p)=-(x+p-1)\Gamma_p(x+p-1)=(-1)^2(x+p-1)(x+p-2)\Gamma_p(x+p-2)=\ldots
\]
\[ \ldots= (-1)^{p-1}(x+p-1)(x+p-2)\ldots(x+1)\Gamma_p(x+1)=(-1)^{p}(x+p-1)(x+p-2)\ldots(x+1)\Gamma_p(x).\]
Consequently 
\begin{equation}
 \forall x \in p\mathbb{Z}_p \quad g(x+p)=f(x)g(x),
 \label{eqpZp}
\end{equation}
where $f(x)=(-1)^{p}(x+p-1)(x+p-2)\ldots(x+1)$.
We compute
\[
\begin{array}{rl}
 g'(x+p)&=f'(x)g(x)+f(x)g'(x),\\
 g''(x+p)&=f''(x)g(x)+2 f'(x)g'(x)+ f(x)g''(x),\\
 &\ldots\\
 g^{(n)}(x+p)&=\sum_{k=0}^n f^{(n-k)}(x)g^{(k)}(x).
\end{array}
\]
By (\ref{eqP}) and (\ref{eqpZp}) for every $x \in p \mathbb{Z}_p$ we have 
\begin{equation}
 P\Big(x+p, f(x)g(x),  f'(x)g(x)+f(x)g'(x), \ldots, \sum_{k=0}^n f^{(n-k)}(x)g^{(k)}(x)\Big)=0. 
\end{equation}
We define $Q(X,Y_0,Y_1, \ldots, Y_n) :=P\big(x+p, f(x)Y_0,  f'(x)Y_0+f(x)Y_1, \ldots, \sum_{k=0}^n f^{(n-k)}(x) Y_k\big) $. Then
$Q(x, g(x),  g'(x), \ldots,  g^{(n)}(x))=0$ for every $x \in p\mathbb{Z}_p$. Since $P$ is minimal, $P$ must divide $Q$. Applying the euclidean algorithm to $LT(P)$ and $LT(Q) = q(X+p)\cdot f(X)^{a_0+a_1+\ldots +a_n}Y_0^{a_0}Y_1^{a_1}\ldots Y_n^{a_n}$ we conclude that there exists $R(X) \in \mathbb{C}_p[X]$ such that
\begin{equation}
 Q(X,Y_0,Y_1, \ldots, Y_n)=R(X)P(X,Y_0,Y_1, \ldots, Y_n) 
 \label{qrp}
\end{equation}
More precisely $R(X)=\frac{q(X+p)}{q(X)} \cdot f(X)^{a_0+a_1+\ldots +a_n}$. Since $\deg f \geq 1$ and $a_0+a_1+\ldots +a_n \neq 0$, then $\deg R \geq 1$.
Since $\mathbb{C}_p$ is algebraically closed, then a nonzero polynomial $R$ must have a root $x_0 \in \mathbb{C}_p$.
Substituting $x_0$ into (\ref{qrp}) we obtain
\[P\big(x_0+p, f(x_0)Y_0,  f'(x_0)Y_0+f(x_0)Y_1, \ldots, \sum_{k=0}^n f^{(n-k)}(x_0) Y_k\big)=0.\]
If $f(x_0)\neq 0$, then $X-x_0-p$ divides $P$, which contradicts the minimality of $P$. Therefore $f(x_0)=0$.
A change of variables then yields 
\[ P(x_0+p,0, Z_1, \ldots, Z_n)=0,\]
for some $x_0 \in \{1-p, 2-p, \ldots, -1\}$.
Substituting $X=x_0+p$ and $Y_0=0$ into (\ref{qrp}) gives us
\[ P\Big(x_0+p, 0, f(x_0+p)Y_1,  \ldots, \sum_{k=1}^n f^{(n-k)}(x_0+p) Y_k\Big) = R(x_0+p)P(x_0+p,0, Y_1, \ldots, Y_n)=0.\]
After performing a suitable change of variables we obtain
\[ P(x_0+2p,0, Z_1, \ldots, Z_n)=0.\]
By induction
\[ \forall m \in \mathbb{Z}, m \geq 1 \quad P(x_0+mp,0, Z_1, \ldots, Z_n)=0.\]
Hence $P(X,0, Z_1, \ldots, Z_n)=0$ and $Y$ divides $P$. We obtain a contradiction with the minimality of $P$. Thus $g$ is differentially transcendental over $\mathbb{C}_p(X)$. As a result $\Gamma_p$ is differentially transcendental over $\mathbb{C}_p(X)$.
$\Box$\\

 \end{document}